\documentclass[12pt]{amsart}
\usepackage{graphicx}
\usepackage{amsmath,amsthm, amsfonts,amssymb, stmaryrd,yfonts,pxfonts,pifont,eufrak, bbm}
\newtheorem{Cor}{Corollary}
 \newtheorem{Lemma}{Lemma}
 
 \newtheorem{ex}{Example}
 \newtheorem{Proposition}{Proposition}
 \theoremstyle{definition}
 
 \theoremstyle{remark}
 \newtheorem{Remark}[Lemma]{Remark}
 \numberwithin{equation}{subsection}

\begin{document}
\title[ON GRAPH INDUCED SYMBOLIC SYSTEMS]{ON GRAPH INDUCED SYMBOLIC SYSTEMS}%
\author{Prashant Kumar and Puneet Sharma}
\address{Department of Mathematics, I.I.T. Jodhpur, NH 65, Nagaur Road, Karwar, Jodhpur-342037, INDIA}%
\email{puneet@iitj.ac.in, kumar.48@iitj.ac.in}%




\begin{abstract}
In this paper, we investigate a shift arising from graph $G$. We prove that any $k$-dimensional shift of finite type can be generated through a $k$-dimensional graph. We investigate the structure of the shift space using the generating matrices for the shift space. We prove that a two dimensional shift space has a horizontally (vertically) periodic point if and only if it possesses a $(m,n)$-periodic point (for some $m,n\in \mathbb{Z}\setminus \{0\}$). We prove that a shift space is finite if and only if it can be generated by permutation matrices. We study the non-emptiness problem and existence of periodic points in terms of the generating matrices.
\end{abstract}
\maketitle

\section{INTRODUCTION}

Symbolic dynamics originated as a tool to investigate various natural and physical phenomena around us. The convenience of symbolic representation and easier computability of the system has attracted attention of several researchers around the globe and the topic has found applications in various branches of sciences and engineering. In particular, the area has found applications in areas like data storage, data transmission and communication systems to name a few \cite{bruce,shanon,lind1}. The structure and dynamics of a symbolic system can be used to investigate the dynamics of a general dynamical system. In fact, it is known that every discrete dynamical system can be embodied in a symbolic dynamical system (with appropriate number of symbols) \cite{fu}. Consequently, it is sufficient to study the shift spaces and its subsystems to investigate the dynamics of a general discrete dynamical system.

Let $A = \{a_i : i \in I\}$ be a finite set and let $d$ be a positive integer. Let the set $A$ be equipped with the discrete metric and let $A^{\mathbb{Z}^d}$, the collection of all functions $c : \mathbb{Z}^d \rightarrow A$ be equipped with the product topology. Any such function $c$ is called a configuration over $A$. Any configuration $c$ is called periodic if there exists $u\in\mathbb{Z}^d~~(u \neq 0)$ such that $c(v+u)=c(v)~~\forall v\in\mathbb{Z}^d$. The set $\Gamma_c= \{w\in {\mathbb{Z}}^{d} : c(v+w)=c(v)~~\forall v\in\mathbb{Z}^d\}$ is called the lattice of periods for the configuration $c$. The function $\mathcal{D} : A^{\mathbb{Z}^d} \times A^{\mathbb{Z}^d} \rightarrow \mathbb{R}^+$ be defined as $\mathcal{D} (x,y) = \frac{1}{n+1}$, where $n$ is the least non-negative integer such that $x \neq y$ in $R_n = [-n,n]^d$, is a metric on $A^{\mathbb{Z}^d}$ and generates the product topology. For any $a\in \mathbb{Z}^d$, the map $\sigma_a : A^{\mathbb{Z}^d} \rightarrow A^{\mathbb{Z}^d}$ defined as $(\sigma_a (x))(k)= x(k+a)$ is a $d$-dimensional shift and is a homeomorphism. For any $a,b\in \mathbb{Z}^d$, $\sigma_a \circ \sigma_b = \sigma_b \circ \sigma_a$ and hence $\mathbb{Z}^d$ acts on $A^{\mathbb{Z}^d}$ through commuting homeomorphisms. For any nonempty $S\subset \mathbb{Z}^d$, any element of $A^S$ is called a pattern over $S$. A pattern is said to be finite if it is defined over a finite subset of $\mathbb{Z}^d$. A pattern $q$ over $S$ is said to be extension of the pattern $p$ over $T$ if $T\subset S$ and $q|_T=p$. The extension $q$ is said to be proper extension if $T\cap Bd(S)=\phi$, where $Bd(S)$ denotes the boundary of $S$. It may be noted that any $k$- dimensional pattern can be visualized as an adjacent placement of some $k-1$- dimensional patterns. For $k-1$-dimensional pattens $B_1,B_2,\ldots,B_r$, let $B=[B_1 B_2 \ldots B_r]_i$ denote the $k$-dimensional pattern obtained by placing $B_1,B_2,\ldots,B_r$ adjacently in the $i$-th direction. We say that a patten $C=[C_1 C_2\ldots C_r]_i$ overlaps progressively with $B=[B_1 B_2 \ldots B_r]_i$ in the $i$-th direction if $B_2B_3\ldots B_r=C_1C_2\ldots C_{r-1}$. Let $\mathcal{F}$ be a given set of finite patterns (possibly over different subsets of $\mathbb{Z}^d$) and let $X=\overline{\{x\in A^{\mathbb{Z}^d}: \text{any pattern from~~} \mathcal{F} \text{~~does not appear in~~} x \}}$. The set $X$ defines a subshift of $\mathbb{Z}^d$ generated by set of forbidden patterns $\mathcal{F}$. If the shift space $X$ can be generated by a finite set of finite patterns, we say that the shift space $X$ is a shift of finite type. We say that a pattern is allowed if it is not an extension of any forbidden pattern. We denote the shift space generated by the set of forbidden patterns $\mathcal{F}$ by $X_{\mathcal{F}}$. Two forbidden sets $\mathcal{F}_1$ and $\mathcal{F}_2$ are said to be equivalent if they generate the same shift space, i.e. $X_{\mathcal{F}_1}= X_{\mathcal{F}_2}$. Refer \cite{bruce, lind1} for details.\\

Let $X$ be a two dimensional shift space over alphabet $\mathcal{A}$ and let $\mathcal{B}_{(M,N)}(X)$ denote the collection of all $M \times N $ patterns allowed for the shift space $X$. Then, $\beta_{(M,N)}: X \rightarrow  (\mathcal{B}_{(M,N)}(X))^{\mathbb{Z}^2}$ defined as $(\beta_{[(M,N)]}(x))_{[(i,j)]}= x_{[i,i+M-1] \times [j,j+N-1]}$ is called  $(M,N)$-higher block code. It can be proved that $\beta_{(M,N)}(X)$ is a shift space (Proposition \ref{hbc}).  Further, it may be noted that for any configuration $c$ in the shift space $X$, any rectangular patterns of size  $M \times N$ appearing in $\beta_{(M,N)}(c)$ placed adjacently (in any direction) overlap progressively (in that direction). A two dimensional shift space of finite type $X_{\mathcal{F}}$ is said to be $(m,n)$-step shift if it can be described by a forbidden set consisting of rectangles of size $(m+1) \times (n+1)$. If the shift space can be described by a forbidden set consisting of blocks of size $1 \times (m+1)$ or $(m+1)\times 1$ , then the shift space $X_{{\mathcal{F}}}$ is called a $m$-step shift. Analogously, for $P=(P_1,P_2,\ldots, P_k)\in\mathbb{N}^k$,  one can define $\mathcal{B}_{P}(X)$ denote the collection of all $ P_1\times P_2\times \ldots \times P_k$ patterns allowed for a $k$-dimensional shift space $X$. Then, Then, $\beta_{P}: X \rightarrow  (\mathcal{B}_{P}(X))^{\mathbb{Z}^k}$ defined as $(\beta_{P}(x))_{[(i_1,i_2,\ldots,i_k)]}= x_{[i_1,i_1+P_1-1] \times [i_2,i_2+P_2-1]\times\ldots \times [i_k,i_k+P_k-1]}$ is called  $(P_1,P_2,\ldots, P_k)$-higher block code(or $P$-higher block code). One again, it can be proved that $\beta_{P}(X)$ is a shift space (Corollary \ref{cc}) and the results (observations) made for the two dimensional case extend analogously for a $k$-dimensional shift space.  A shift space $X_{\mathcal{F}} $ is said to be aperiodic if it does not contain any periodic points.\\ 

Let $G$ be a graph with finite set vertices $V$ and finite set of edges $E$. It can be seen that the set of bi-infinite walks over a graph is a $1$-step shift of finite type. Also, for any given shift of finite type $X$, there exists a higher block shift (conjugate to $X$) which can be generated by a finite graph $G$. Consequently, every one dimensional shift of finite type can be visualized as a shift generated from some graph \cite{bruce,lind1}.\\

For multidimensional shifts of finite type, it is known that given a set of forbidden patterns, the non-emptiness problem for multidimensional shift spaces is undecidable \cite{ber}. In \cite{emma}, the authors show that the sets of periods of multidimensional shifts of finite type are exactly the sets of integers of the complexity class NE. They also give characterizations for general sofic and effective subshifts. In \cite{coven}, authors prove that a multidimensional of finite type has a power that can be realized as the same power of a tiling system. They show that the set of entropies of tiling systems equals the set of entropies of shifts of finite type. It is known that multidimensional shifts of finite type with positive topological entropy cannot be minimal\cite{quas}. Infact, if $X$ is subshift of finite type with positive topological entropy, then $X$ contains a subshift which is not of finite type, and hence contains infinitely many subshifts of finite type \cite{quas}. In \cite{hoch4}, Hochman proved that $h\geq 0$ is the entropy of a $\mathbb{Z}^d$ effective dynamical system if and only if it is the lim inf of a recursive sequence of rational numbers. For two dimensional shifts, Lightwood proved that strongly irreducible shifts of finite type have dense set of periodic points \cite{sam}.  In \cite{pd}, the authors characterized a multidimensional shift of finite type using an infinite matrix. In \cite{pd1}, authors gave an algorithmic approach to address the non-emptiness problem for multidimensional shift space. They give an algorithm to generate the elements of the shift space using finite matrices. In the process, they prove that that elements of d-dimensional shift of finite type can be characterized by a sequence of finite matrices.\\

Although a lot of work for multidimensional shift spaces has been done, graph induced multidimensional shifts have not been investigated. If $\{G_1,G_2,\ldots,G_d\}$ is a set of $d$ graphs with a common set of vertices $V$, the collection naturally induces a $d$-dimensional shift of finite type (where $i$-th graph determines the compatibility of the vertices in the $i$-th direction). In this paper, we investigate the relation between the structure of the generating graphs $G_i$ and the shift space generated. In particular, we answer some of the questions relating the the structure of the underlying graphs with the non-emptiness problem of the shift space and existence of periodic points. For example, can every shift of finite type $X$ be generated by a finite set of graphs? when does a given collection $\{G_1,G_2,\ldots,G_d\}$ of graphs generate a non-empty shift space? When does a multidimensional shift generated by $\{G_1,G_2,\ldots,G_d\}$ exhibit periodic points? Does existence of periodicity in one direction ensure the periodicity in other directions? We now give answers to some of these questions relating the multidimensional shift space and the generating set of graphs.

\section{Main Results}

\begin{Proposition} \label{odg}
For any two dimensional one step shift of finite type $X$, there exists a two dimensional graph $G$ such that $X=X_{G}$.	
\end{Proposition}

\begin{proof}
Let $X$ be a two dimensional one step shift of finite type over the finite alphabet set $\mathcal{A}$. As $X$ is one step, $X$ is generated by a forbidden set $\mathcal{F}$ such that any element of $\mathcal{F}$ is of the form $^{b}_{a}$ or $ab$ (where $a,b\in \mathcal{A}$). Define a graph $H$ ($V$) with $\mathcal{A}$ as the set of vertices and $\exists $ a directed edge from vertex $a$ to vertex $b$ in $H$ ($V$) if and only if $ab$ ($^{b}_{a}$) does not belong to $\mathcal{F}$. Then, as $G=(H,V)$ is a two dimensional graph that captures horizontal and vertical compatibility of the elements of $\mathcal{A}$, $G$ generates any arbitrary element of $X$. Consequently, $X=X_G$ and the proof is complete.
\end{proof}	

\begin{Remark}
The above result establishes that any two dimensional one step shift of finite type can be generated by a two dimensional graph. It may be noted that for a $k$-dimensional one step shift of finite type $X$, if $H_i$ is the graph that captures the compatibility of the symbols in the $i$-th direction, then similar arguments establish that $G=(H_1,H_2,\ldots,H_k)$ generates an arbitrary element of $X$ (and conversely). Consequently, the above result holds for any higher dimensional one step shift and we get the following corollary.
\end{Remark}

\begin{Cor}\label{kdg}
	For any $k$-dimensional one step shift of finite type $X$, there exists a $k$- dimensional graph $G$ such that $X=X_{G}$.	
\end{Cor}

\begin{proof}
The proof follows from discussions in Remark 1. 	
\end{proof}	

\begin{Proposition} \label{hbc}
For any two dimensional shift space $X_{\mathcal{F}}$, $X^{(M,N)} $ is a shift space conjugate to $X_{\mathcal{F}}$.
\end{Proposition}

\begin{proof}
Let $X_{\mathcal{F}}$ be a shift space generated by the forbidden set $\mathcal{F}$ and let ${(M,N)} \in \mathbb{N}^{2} $. Let $\mathcal{F}^*$ be the set obtained by replacing any forbidden pattern $P$ of size less than size $M\times N$ by all $M\times N$ extensions of $P$. Then, $X_{\mathcal{F}}= X_{\mathcal{F}^*}$ and hence we obtain a modified forbidden set generating $X_{\mathcal{F}}$ such that all the forbidden patterns in the generating forbidden set are bigger than a rectangle of size $M\times N$. Further, as all the forbidden patterns can be extended to rectangles of uniform size to generate the same space, we assume all the elements of the forbidden set to be rectangles of size $R \times S$ (for some integers $R,S\in\mathbb{N}$).

For any $P\in F$, define $P^{(M,N)}$ to be a pattern of size $(R-M+1)\times (S-N+1)$ over $(\mathcal{B}_{(M,N)}(X))$ defined as $P^{(M,N)}_{[(k,l)]}= P_{[k,k+M-1]\times [l,l+N-1]}$, i.e. the $M\times N$ rectangle with left bottom corner at $(k,l)$ is placed at $(k,l)$. Let $\mathcal{F}_{1}=\{P^{(M,N)}: P\in \mathcal{F}\}$. Further, let  $\mathcal{F}_{2}= \{P_{1}P_{2}:P_{1},P_{2} \in \mathcal{A}_{X}^{[M,N]} \ such  \ that \ P_{1}\ and\ P_{2}  \ do \ not \ overlap \ progressively \ horizontally\}$ and let   $\mathcal{F}_{3}= \{^{P_{2}}_{P_{1}} \ : \ P_{1},P_{2}\in \mathcal{A}_{X}^{[M,N]} \ such  \  that \ P_{1}\ and \ P_{2}  \ do \ not \ overlap \ progressively \\ vertically \}$.

Note that as elements of $\mathcal{F}$ are forbidden for $X$, elements of $\mathcal{F}_{1}$ are forbidden for $X^{(M,N)}$. Also, as any two blocks placed adjacently for $X^{(M,N)}$ must overlap progressively, $\mathcal{F}_{2}$ and $\mathcal{F}_{3}$ are also forbidden for $X^{(M,N)}$ and thus $X^{(M,N)} \subset \cap_{i=1}^3 X_{\mathcal{F}_{i}}$ or $ X^{(M,N)} \subseteq X_{{\mathcal{F}_{1}} \cup {\mathcal{F}_{2}}  \cup {\mathcal{F}_{3}}  } $. Conversely, for any element $x$ in $X_{{\mathcal{F}_{1}}\cup {\mathcal{F}_{2}}\cup {\mathcal{F}_{3}}  }$, as adjacent placement of blocks not overlapping progressively is forbidden, any two adjacent blocks overlap progressively. Further, as elements of $\mathcal{F}_1$ are forbidden for $X_{{\mathcal{F}_{1}} \cup {\mathcal{F}_{2}}  \cup {\mathcal{F}_{3}}}$, any block forbidden for $X$ does not appear in $x$ .  Consequently, $x\in X^{(M,N)} $ and the proof for $X^{(M,N)}= X_{{\mathcal{F}_{1}}\cup {\mathcal{F}_{2}}\cup {\mathcal{F}_{3}}}$ is complete.

Further,for any $x\in X$ as $((\beta_{(M,N)})(x))_{(i,j)}$ is a $M\times N$ pattern of $x$ with left corner at $x_{(i,j)}$,  the map $\beta_{(M,N)}$ defines a conjugacy between shift space X and $ X^{(M,N)}$.
\end{proof}

\begin{Remark}\label{kde}
The above result establishes that any two dimensional shift is conjugate to its higher block code $X^{(M,N)}$. The proof uses the fact slicing any given configuration in patterns of size $M\times N$ at each $(r,s)\in \mathbb{Z}^2$ (and placing it at each $(r,s)\in \mathbb{Z}^2$) yields an element of $(\mathcal{B}_{(M,N)}(X))^{\mathbb{Z}^2}$. The correspondence is natural and indeed is a conjugacy between $X$ and $X^{(M,N)}$. Note that if $X$ is a $k$-dimensional shift space and $P\in\mathbb{N}^k$, then slicing any configuration in $X$ in patterns of size $P$ at each point in $\mathbb{Z}^k$ (and placing the slice at each point in $\mathbb{Z}^k$) extends the above result for a $k$- dimensional shift space. Thus we get the following corollary.
\end{Remark}

\begin{Cor} \label{cc}
For any $k$-dimensional shift space $X_{\mathcal{F}}$ and $P\in\mathbb{N}^k$, $X^{P} $ is a shift space conjugate to $X_{\mathcal{F}}$.
\end{Cor}

\begin{proof}
The proof follows from discussions in Remark \ref{kde}.
\end{proof}

\begin{Proposition}
For any two dimensional shift space of finite type $X_{\mathcal{F}}$, there exists a graph $G$ such that $X = X_{G}$.
\end{Proposition}

\begin{proof}
Let $X_{\mathcal{F}}$ be a shift space of finite type generated by the forbidden set $\mathcal{F}$. If all the elements of $\mathcal{F}$ are of type $\{\alpha \beta\} $ or $\{^{\alpha}_{\beta}\}$, then $X_{\mathcal{F}}$ is one step shift of finite type. If not, let all the elements of $\mathcal{F}$ be rectangles of size $M\times N$.
By previous proposition,  Since $ X^{[M,N]}$ can be viewed as one step of finite type over alphabet $\mathcal{A}_{X}^{[(M,N)]}$, shift space $X_{\mathcal{F}}$ can be expressed as one step shift of finite type. For $\mathcal{V}= \mathcal{B}_{(M,N)}(X)$, define the graph $H_1=(\mathcal{V},E_1)$ as a graph with set of vertices $\mathcal{V}$ where any two elements of $\mathcal{V}$ are connected if they overlap progressively horizontally. Let $H_2=(\mathcal{V},E_2)$ be the graph with $\mathcal{V}$ as the set of vertices where any two elements of $\mathcal{V}$ are connected if they overlap progressively vertically. Then $G=(H_1,H_2)$ generates $X^{(M,N)}$ and the proof is complete.
\end{proof}

\begin{Remark} \label{kdc}
The above result establishes that any two dimensional shift of finite type can be generated by a two dimensional graph. The proof uses the fact as any shift of finite type is conjugate to $X^{(M,N)}$, the shift space can be visualized as a one step shift of finite type. Further, as any one step shift of finite type can be generated through a graph, any two dimensional shift of finite type can be generated by some graph $G=(H,V)$. As any $k$-dimensional one step shift of finite type can be generated through a graph (Corollary \ref{kdg}), any $k$-dimensional shift of finite type can be generated by a $k$-dimensional graph. Consequently, an analogous extension of the above result is true and we get the following corollary.   
\end{Remark}

\begin{Cor}
Every $k$-dimensional shift of finite type $X_{\mathcal{F}}$ can generated by some $k$-dimensional graph $G$.
\end{Cor}

\begin{proof}
The proof follows from discussions in Remark \ref{kdc}.
\end{proof}

\begin{ex}
Let X be two dimensional shift space with alphabet $\{e,f,g\}$ with forbidden pattern set $\mathcal{F}= \{ff,gg,fe,eg, \ ^{f}_{f}, \ ^{e}_{e}, \ ^{g}_{g}, \ ^{e}_{f} , \  ^{g}_{e} \}$. Then, graph G for this shift space is given by figure 1.
\begin{figure}[h]
		\includegraphics[height=5.0cm, width=10.0cm]{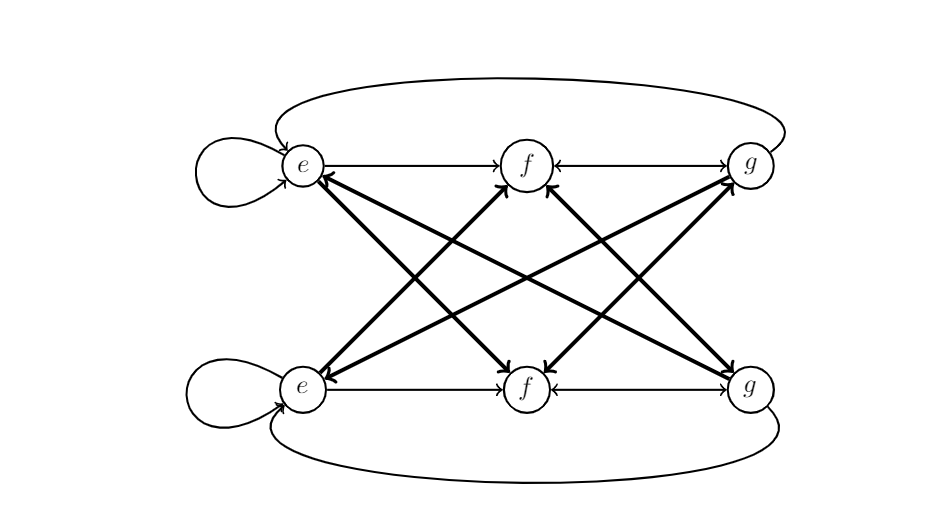}
		\caption{}
\end{figure}

Then, as there exists $2\times 2$ patterns whose infinite repetition (in both directions) tiles the plane in an allowed manner, the shift space is non-empty and exhibits periodic points. Further, an arbitrarily large central block of a given configuration can be infinitely repeated to obtain an element of $X$, the shift space exhibits a dense set of periodic points.

\end{ex}

\begin{ex}
	Let X be two dimensional Golden Mean shift space over alphabet $\{0,1\}$ with forbidden pattern set $\mathcal{F}= \{11, \ ^{1}_{1}  \}$. Then, $X= X_{G}$, where graph G is given by figure 2.
	\begin{figure}[h]
		\includegraphics[height=5.0cm, width=10.0cm]{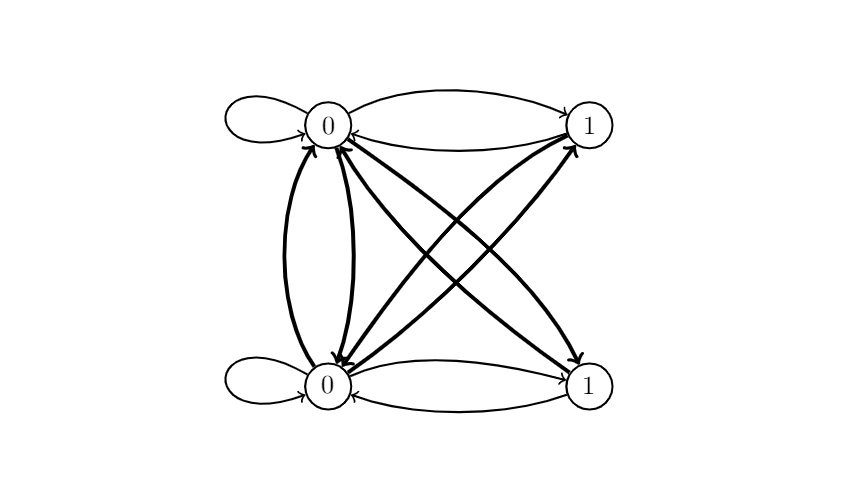}
		\caption{}
	\end{figure}

Then, appearance of two consecutive $1$'s is forbidden in any direction. As the configuration comprising of all $0$'s is a valid element of $X$, the shift space $X$ is indeed non-empty. Note that any allowed $2\times 2$ pattern can be extended to a valid repetition of the shift space $X$. Once again, as arbitrarily large central blocks of a given configuration can be infinitely repeated to obtain an element of $X$, the shift space exhibits a dense set of periodic points.
\end{ex}


\begin{Proposition}
For any one step $2$-dimensional shift of finite type $X$, $X$ has a horizontally periodic point if and only if $X$ has a $(m,n)$ periodic point (for some $m,n\in \mathbb{Z}\setminus \{0\}$).
\end{Proposition}

\begin{proof}
Let $X$ be a one step shift of finite type and let $x\in X$ be a $(m,0)$ periodic point. Then, note that $x$ is a infinite horizontal repetition of an infinite vertical strip of width $m$ (say $\mathbb{S}$). Further as $\mathbb{S}$ can be realized as a vertical arrangement of one dimensional strips of length $m$, there exists a $1\times m$ block $a_1,a_2\ldots a_m$ which appears twice in $\mathbb{S}$ (say at heights $u$ and $v$). Consequently, infinite repetition of the block $x_{[0,m-1]\times [u,v-1]}$ is an element of $X$ and is periodic of period $(m,v-u)$. \\

Conversely, if $X$ has a $(m,n)$ periodic point then there exists an infinite (horizontal) strip $\mathbb{S}$ such that $x$ is a vertical arrangement of shifts of $\mathbb{S}$ (where $\sigma^{(-m,0)}(\mathbb{S}), \mathbb{S}, \sigma^{(m,0)}(\mathbb{S}), \sigma^{(2m,0)}(\mathbb{S}), \ldots$ are placed vertically one over the other to obtain $x$). As the blocks of size $m\times n$ are finite, there exists a block $B_0$ of size $m\times n$ that appears in $x$  at  $(u,0)$ and $(v,0)$. Consequently, if $B_0 B_1\ldots B_k B_0$ is a block appearing in $X$ then the $k \times k$ rectangular arrangement of $B_0,B_1,\ldots,B_k$ where $B_{(k-j+i+1) \text{mod}(k+1)}$ is placed at $(i,j)$-th position is an allowed rectangular block. Further, as infinite repetition of the block generated yields an allowed configuration of $X$, the shift space exhibits a horizontally periodic point and the proof is complete.
\end{proof}

\begin{Remark}\label{pp}
The above proof establishes equivalence of existence of periodic points  with existence of horizontally periodic points for a shift of finite type. The proof uses the fact that any $(m,n)$ periodic point (with $m,n\neq 0$) can be realized as a vertical arrangement of shifts of an infinite horizontal strip of height $n$. The periodic point generated is also vertically periodic and hence the proof establishes equivalence of existence of periodic points  with existence of vertically periodic points for a shift of finite type. Thus we get the following corollary.
\end{Remark}

\begin{Cor}
For any one step $2$-dimensional shift of finite type $X$, $X$ has a vertically periodic point if and only if $X$ has a $(m,n)$ periodic point (for some $m,n\in \mathbb{Z}\setminus \{0\}$).
\end{Cor}

\begin{proof}
The proof follows from discussions in  Remark \ref{pp}.
\end{proof}

%

%
%
%
 %
	
\begin{Proposition}
	A two dimensional shift space $X_{G}$ is finite iff it can be generated by a pair of permutation matrices.
\end{Proposition}

\begin{proof}
Firstly note that any finite shift space is a union of finitely many periodic points (with finite orbits). Also, if $X$ itself is a single periodic orbit then $X$ can be visualized as an infinite repetition (both horizontal and vertical) of an $m\times n$ rectangle. Then, if $H$ and $V$ are indexed with allowed rectangles of size $m\times n$ capturing horizontal and vertical compatibility of the indices then $H$ and $V$ are permutation matrices and the graph $G=(H,V)$ generates the shift space $X$. Finally if $X$ is a union of periodic orbits, a similar argument applied to each periodic orbit (and collating the set of indices to generate $H$ and $V$) generates a pair of permutation matrices that generate $X$ and the proof of forward part is complete.

Conversely, if the generating matrices are permutation matrices then fixing the entry at the origin fixes the entries in the immediate neighborhood and hence fixes all the entries at other coordinates. Consequently, the shift space $X$ is finite and the proof is complete.
\end{proof}	

\begin{Remark}
The above result establishes that a two dimensional shift space is finite if and only if it can be generated by a pair of permutation matrices. However, finiteness of the shift space $X_{G}$ does not enforce the generating matrices $H$ and $V$ to be permutation matrices. To establish our claim, let $X$ be the shift space generated by the graph shown in Figure $3$. Then, it can be seen that although the shift generated by the graph is finite, the associated adjacency matrices H and V are not permutation matrices and hence the claim is indeed true.

\begin{figure}[h]
	\includegraphics[height=6.0cm, width=12.0cm]{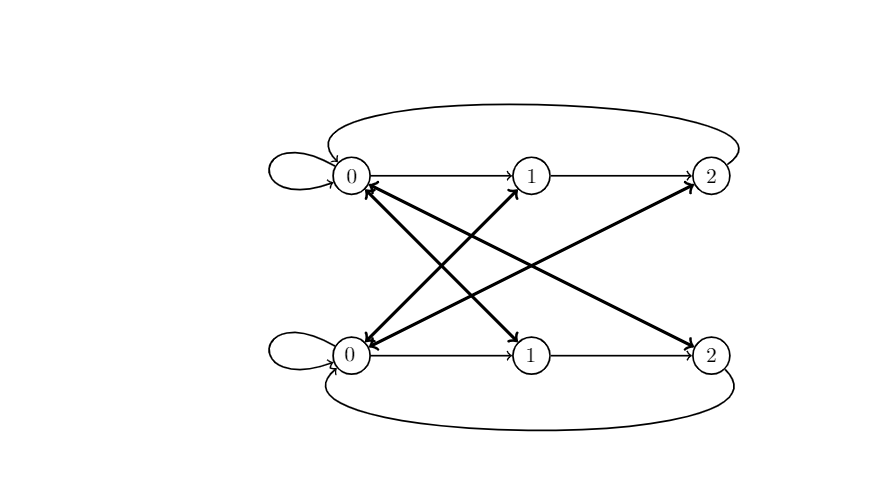}
	\caption{}
\end{figure}

$$\textit{H}= \bordermatrix{ & 0 & 1 & 2 \cr
	0 & 1 & 1 & 0 \cr
	1 & 0 & 0 & 1 \cr
	2 & 1 & 0 & 0 \cr 	
}
\ \ \ \ \ \ \ \ \
\textit{V}= \bordermatrix{ & 0 & 1 & 2 \cr
	0 & 0 & 1 & 1 \cr
	1 & 1 & 0 & 0 \cr
	2 & 1 & 0 & 0 \cr 	
}$$

But $X_{G}$ is finite as it is the orbit of a single periodic point (given below):

$$ {\begin{array}{ccccccccccccccccccccccc}
	
	\ldots & \vdots & \vdots & \vdots & \vdots & \vdots  & \vdots & \vdots & \vdots   & \vdots & \vdots & \vdots & \vdots & \vdots & \vdots  & \vdots & \vdots & \vdots   & \ldots   \\

	\ldots & 1 & 2 & 0 & 0 & 1  & 2 & 0 & 0  & 1 & 2 & 0 & 0 & 1  & 2 & 0 & 0  & \ldots   \\
	
	\ldots &	0 & 0 & 1 & 2 & 0 & 0 &	1 & 2 &	0 & 0 & 1 & 2 & 0 & 0 &	1 & 2 & \ldots    \\
	
	\ldots & 1 & 2 & 0 & 0 & 1  & 2 & 0 & 0  & 1 & 2 & 0 & 0 & 1  & 2 & 0 & 0  & \ldots   \\
	
	\ldots &	0 & 0 & 1 & 2 & 0 & 0 &	1 & 2 & 	0 & 0 & 1 & 2 & 0 & 0 &	1 & 2 & \ldots    \\
	
	\ldots & 1 & 2 & 0 & 0 & 1  & 2 & 0 & 0  & 1 & 2 & 0 & 0 & 1  & 2 & 0 & 0  & \ldots   \\
	\ldots & \vdots & \vdots & \vdots & \vdots & \vdots  & \vdots & \vdots & \vdots   &\vdots & \vdots & \vdots & \vdots & \vdots & \vdots  & \vdots & \vdots & \vdots   &  \ldots   \\	
	\end{array} } $$

Consequently, finiteness of the shift space $X$ does not guarantee the generating matrices to be permutation matrices.\\
\end{Remark}

We now discuss non-emptiness of shift spaces using adjacency matrices H and V. Note that while $(HV)_{ij}$ computes number of ways pattern $^{\ \ \ {j}}_{ \  {i} \ \ \ }$  can be extended to triangular pattern of form $^{\ \ \ {j}}_{ \  {i} \ {k} } \in \mathcal{B}(X_{G})$, $(VH)_{ij}$ computes the number of ways pattern $^{\ \ \ {j}}_{ \  {i} \ \ \ }$  can be extended to triangular pattern of form $^{\  {l} \ \ {j}}_{ \  {i} \  } \in \mathcal{B}(X_{G})$. As removing the vertices with no incoming (or outgoing) edge (both horizontally or vertically) does not effect the shift space generated, we assume that the generating matrices do not contain any zero row or zero column.

\begin{Proposition}
Let G be graph with associated adjacency matrices H and V. If H and V are irreducible permutation matrices then  $HV = VH$ iff $X_{G} \neq \phi$.
\end{Proposition}
\begin{proof}
Let $X_G$ be the shift space generated by $G=(H,V)$. If $H$ and $V$ are permutation matrices then fixing an entry at origin uniquely determines the immediate neighbors of any symbol (in both horizontal and vertical directions). Further as $HV$ and $VH$ are permutation matrices (characterizing blocks of form ${\begin{array}{cc} &  b  \\ a & *\\ \end{array} }$ and ${\begin{array}{cc} * &  b  \\ a &  \\ \end{array} }$ respectively), any block ${\begin{array}{cc} &  b  \\ a & *\\ \end{array} }$ can be extended to a $2\times 2$ square if and only if $HV=HV$. As the possible choices for immediate neighbor are unique, $X_G$ is non-empty if and only if $HV=VH$ and the proof is complete.
\end{proof}

\begin{Remark}\label{ne}
It may be noted that for any shift generated by permutation matrices, as the immediate neighborhood of a symbol is uniquely determined, the shift space generated by permutation matrices is always finite (may be empty). The above result establishes that a two dimensional shift space generated by a pair of irreducible matrices is non-empty if and only if and only if the generating matrices commute with each other. The proof follows from the fact that if $HV=VH$, any pattern of the form ${\begin{array}{cc} &  c  \\ a & b\\ \end{array} }$ can be extended to a $2\times 2$ square and hence shift space generated is non-empty ( in fact, is a finite shift space comprising of a single periodic orbit). Note that if $(HV)_{ij} \neq 0 \Leftrightarrow (VH)_{ij} \neq 0~~\forall i,j$ then the shift space is non-empty and hence a more general form of the above result is true.  Further, note that if $(HV)_{ij}\neq 0\implies (VH)_{ij}\neq 0$ (or $(VH)_{ij}\neq 0\implies (HV)_{ij}\neq 0$), shift space generated does not contain any forbidden pattern of the form ${\begin{array}{cc} &  c  \\ a & b\\ \end{array} }$ (or ${\begin{array}{cc} a &  b  \\ c &  \\ \end{array} }$) and hence the shift space generated is once again non-empty. Thus we get the following corollaries.
\end{Remark}

\begin{Cor}
	Let $X_{G}$ be a shift space generated by $G=(H,V)$. If $(HV)_{ij} \neq 0 \iff (VH)_{ij} \neq 0~~\forall i,j$ then $X_{G} \neq \phi$.
\end{Cor}

\begin{proof}
	The proof follows from the discussions in Remark \ref{ne}.
\end{proof}

\begin{Cor}
	Let $X_{G}$ be a shift space  generated by $G=(H,V)$. If $(HV)_{ij} \neq 0 \Rightarrow (VH)_{ij} \neq 0$. Then $X_{G} \neq \phi$.
\end{Cor}

\begin{proof}
The proof follows from the fact that if $(HV)_{ij} \neq 0 \Rightarrow (VH)_{ij} \neq 0$, the shift space does not contain forbidden pattern of the form  ${\begin{array}{cc} &  c  \\ a & b\\ \end{array} }$. Consequently, any arbitrarily large $1\times r$ pattern can be extended to an $r\times r$ square. As the shift space contains valid arbitrarily large squares, the shift space is non-empty and the proof is complete.
\end{proof}

\begin{Remark}
The above proposition proves that if $H$ and $V$ are irreducible permutation matrices, the shift space is non-empty if $HV=VH$. Further as $HV=VH$ ensures that every pattern of the form ${\begin{array}{cc} &  c  \\ a & b\\ \end{array} }$ (or ${\begin{array}{cc} a &  b  \\ c &  \\ \end{array} }$) extendable to $2\times 2$ square, the shift space is non-empty if $HV=VH$. However, if $H$ and $V$ are not irreducible, any element of the shift space can possibly be generated by sub-matrices of $H$ and $V$ and hence the shift space can be non-empty even when $HV=VH$ does not hold good. We now give an example in support of our claim.
\end{Remark}

\begin{ex}
Let $X$ be a shift space generated by the graph in Figure 4.

\begin{figure}[h]
	\includegraphics[height=5.0cm, width=10.0cm]{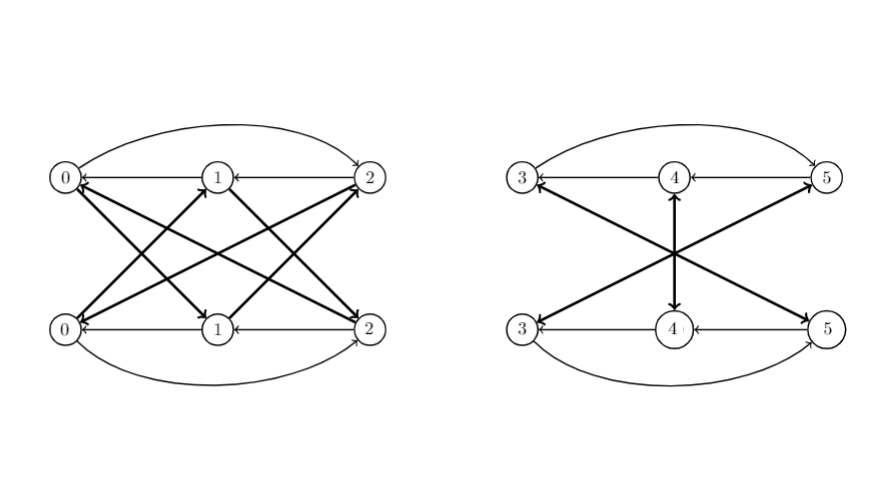}
	\caption{}
\end{figure}

Then,
$$\textit{H}= \bordermatrix{ & 0 & 1 & 2 & 3 & 4 & 5 \cr
	0 & 0 & 0 & 1 & 0 & 0 & 0\cr
	1 & 1 & 0 & 0 & 0 & 0 & 0 \cr
	2 & 0 & 1 & 0 & 0 & 0 & 0\cr
	3 & 0 & 0 & 0 & 0 & 0 & 1\cr
	4 & 0 & 0 & 0 & 1 & 0 & 0\cr
	5 & 0 & 0 & 0 & 0 & 1 & 0\cr
}
\ \ \ \ \ \ \ \ \
\textit{V}= \bordermatrix{ & 0 & 1 & 2 & 3 & 4 & 5 \cr
	0 & 0 & 1 & 0 & 0 & 0 & 0\cr
	1 & 0 & 0 & 1 & 0 & 0 & 0\cr
	2 & 1 & 0 & 0 & 0 & 0 & 0\cr
	3 & 0 & 0 & 0 & 0 & 0 & 1\cr
	4 & 0 & 0 & 0 & 0 & 1 & 0\cr
	5 & 0 & 0 & 0 & 1 & 0 & 0\cr
}$$

Note that $G$ can be written as a union of disjoint graphs $G_1$ and $G_2$ indexed by symbols $0,1,2$ and $3,4,5$ respectively. Further, while matrices capturing horizontal and vertical compatibility of $G_1$ commute, matrices capturing horizontal and vertical compatibility of $G_2$ do not commute and hence $X_{G_1}\neq \phi$ but $X_{G_2}=\phi$. Consequently, $X_G=X_{G_1}$ and the shift space is indeed non-empty. Thus, the shift space is generated by a non-commuting pair of permutation matrices.\\
\end{ex}

\begin{Remark}
The above results investigate the non-emptiness of the shift space using the matrices $HV$ and $VH$. However, as $HV^T$ and $V^TH$ characterizes allowed patterns of the form ${\begin{array}{cc} a &  b  \\  & c\\ \end{array} }$ and ${\begin{array}{cc} a &   \\ b & c \\ \end{array} }$ respectively, the non-emptiness problem and existence of periodic points can be investigated using the matrices $HV^T$ and $V^TH$. However, it is worth mentioning that the two conditions are indeed independent and hence can be used independently to investigate the shift space under discussion. We now give an example in support of our claim.
\end{Remark}

\begin{ex}
Let $X$ be the shift space arising from graph in figure-5 over symbol set $ \{1,2,3 \}$.

\begin{figure}[h]
		\includegraphics[height=5.0cm, width=10.0cm]{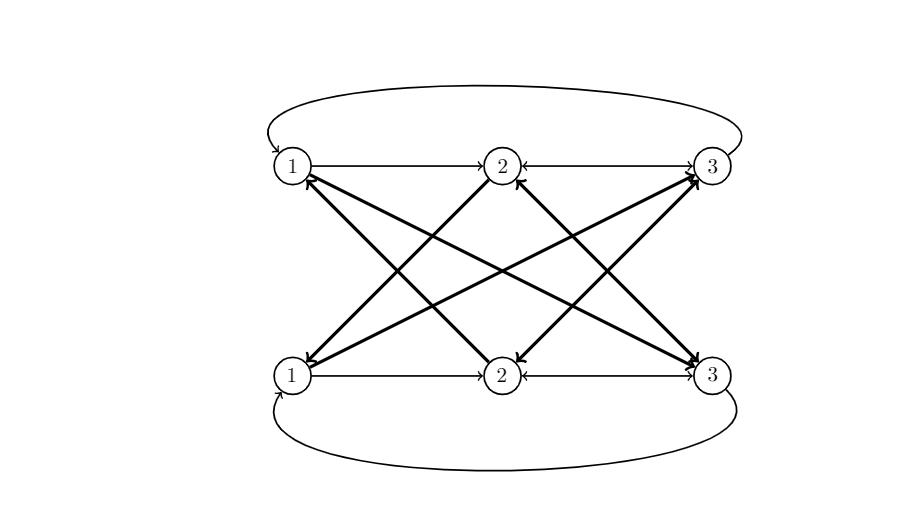}
		\caption{}
	\end{figure}

Then, generating matrices corresponding to given graph are:\\

	$$\textit{H}= \bordermatrix{ & 1 & 2 & 3 \cr
		1 & 0 & 1 & 0 \cr
		2 & 0 & 0 & 1 \cr
		3 & 1 & 1 & 0 \cr 
	}
	\ \ \ \ \ \ \ \ \ 
	\textit{V}= \bordermatrix{ & 1 & 2 & 3 \cr
		1 & 0 & 0 & 1 \cr
		2 & 1 & 0 & 1 \cr
		3 & 0 & 1 & 0 \cr 	
	}$$
	Then,
	$$\textit{HV}= \bordermatrix{ & 1 & 2 & 3 \cr
		1 & 1 & 0 & 1 \cr
		2 & 0 & 1 & 0 \cr
		3 & 1 & 0 & 2 \cr 
	}
	\ \ \ \ \ \ \ \ \ 
	\textit{VH}= \bordermatrix{ & 1 & 2 & 3 \cr
		1 & 1 & 1 & 0 \cr
		2 & 1 & 2 & 0 \cr
		3 & 0 & 0 & 1 \cr 	
	}$$

For the above example, one can find indices $i,j$ such that of $(HV)_{ij} \neq 0$ but $(VH)_{ij} =0 $  (and indices $k,l$ such that $(VH)_{kl} \neq 0 $ but $(HV)_{kl} =0$). Consequently, the condition $(HV)_{ij}=0 $ iff $ (VH)_{ij}=0$ does not hold good and the derived results cannot be used to investigate the non-emptiness of the shift space.   However, as $H=V^T$, we have $HV^{T}=V^{T}H$ and the shift space is indeed non-empty (and possesses periodic points)
\end{ex}

\begin{Proposition}
Let $X_{G}$ be a shift space generated by $G=(H,V)$. If $(HV)_{ij} \neq 0 \iff (VH)_{ij} \neq 0~~\forall i,j$ then $X_{G}$ possesses periodic points (of arbitrarily large periods).
\end{Proposition}

\begin{proof}
Let $X_{G}$ be a shift space generated by $G=(H,V)$ and let $m\in \mathbb{N}$. Let $u$ be a block of size $1\times m$. As $(HV)_{ij} \neq 0 \iff (VH)_{ij} \neq 0~~\forall i,j$, $u$ can extended to a pattern of size $k\times m$ (for any $k\in \mathbb{N}$). Without loss of generality, let $u$ be extended to a rectangle $v$ of size $s\times m$ such that $v_{00}=v_{ms}$. As $v_{00}=v_{ms}$, the block $v$ can be further extended to the block ${\begin{array}{cc}  &  v  \\ v & \\ \end{array} }$  (along the line $sx-my=0$) to obtain a valid pattern of $X$. Further, as $(HV)_{ij} \neq 0 \iff (VH)_{ij} \neq 0~~\forall i,j$, the pattern can be extended to valid $2m\times 2s$ pattern for the shift space. Finally, note that the infinite such repetition of $v$ (along the line $sx-my=0$) and extending the pattern with the same choices (as in the previous step) yields a valid periodic point for the shift space. As the proof holds for any $m\in \mathbb{N}$, the shift space contains periodic points of arbitrarily large periods and the proof is complete.
\end{proof}

\begin{Remark}\label{hvp}
The above result establishes the existence of periodic points under the condition $(HV)_{ij} \neq 0 \iff (VH)_{ij} \neq 0~~\forall i,j$. The proof uses the condition to extend the pattern for the form to a valid $2m\times 2s$ pattern. As such a repetition can be made infinitely often, filling the choices in a unique manner at each step yields a periodic point for the shift space. Note that as such an extension is possible under $(HV)_{ij} \neq 0 \implies (VH)_{ij} \neq 0~~\forall i,j$, the result holds good under a weaker condition. Further, as similar arguments establish the result under the condition $(VH)_{ij} \neq 0 \implies (HV)_{ij} \neq 0~~\forall i,j$, we get the following corollary.
\end{Remark}

\begin{Cor}
Let $X_{G}$ be a shift space generated by $G=(H,V)$. If $(HV)_{ij} \neq 0 \implies (VH)_{ij} \neq 0~~\forall i,j$ (or $(VH)_{ij} \neq 0 \implies (HV)_{ij} \neq 0~~\forall i,j$) then $X_{G}$ possesses periodic points (of arbitrarily large periods).
\end{Cor}

\begin{proof}
The proof follows from discussions in Remark \ref{hvp}.
\end{proof}

\vskip 0.5cm

Let $X$ be a shift space generated by a graph $G$. It may be noted that if $(HV)_{ij}=0$ then any block of the form ${\begin{array}{cc} &  j  \\ i & *\\ \end{array} }$ is forbidden for the shift space $X$. Consequently, the set $\{(i,j): (HV)_{ij}=0 \text{~but~} (VH)_{ij}\neq 0\}$ characterizes all patterns of the form ${\begin{array}{cc} * &  j  \\ i &  \\ \end{array} }$ which cannot be extended to a $2\times 2$ square. Similarly, the set  $\{(i,j): (VH)_{ij}=0 \text{~but~} (HV)_{ij}\neq 0\}$ characterizes all patterns of the form ${\begin{array}{cc} &  j  \\ i & *\\ \end{array} }$ which cannot be extended to a $2\times 2$ square. As such patterns do not contribute towards generation of elements of $X$, ignoring this piece of information leads to generation of elements with reduced complexity. Thus let us set $(HV)_{ij}=0$ if $(VH)_{ij}=0$ (and conversely). \\

Let $A_{1} = \{  ^{\ \ \ {c}}_{ \  {a} \ {b} } \ : \exists~~  d \ \in \mathcal{V}(G) \ such \  that \  ^{d \ \ c}_ { a \ \ b}  \  \in \mathcal{B}({X}_{G}) \}$ and $A_{2} = \{  {^{{y} \ {z}}_{x}    }  \ : \exists \ w  \in \mathcal{V}(G) \ such \  that \  ^{y \ \ z}_ { x \ \ w}  \  \in \mathcal{B}({X}_{G}) \}$. Let $M$ and $N$ be matrices indexed by elements of $A_1$ and $A_2$ in the following manner:\\

For  $I \ = \  ^{\ \ \  \ {a_{3}}}_{ \  {a_{1}} \ {a_{2}} }~~$, $~~J = \  ^{\ \ \  \ {a_{5 }}}_{ \  {a_{3}} \ {a_{4}} }, ~~R =  {^{{b_{2}} \ {b_{3}}}_{b_{1}}} $ and $~~S = {^{{b_{4}} \ {b_{5}}}_{b_{3}}}$\\

$M_{IJ}=
\begin{cases}
	1, & if \ \   {^{{a_{3}} \ {a_{4}}}_{a_{2}}    }  \in {A}_{2}  \\
	
	0, &  otherwise
\end{cases}
~~$
and
$~~N_{RS}=
\begin{cases}
	1, & if \ \   ^{\ \ \  \ {b_{4}}}_{ \  {b_{2}} \ {b_{3}} }  \in {A}_{1}  \\
	
	0, &  otherwise
	
\end{cases}	
$

We identify the pair of indices  $I= ^{\ \ \  \ {a_{3}}}_{ \  {a_{1}} \ {a_{2}} }$ and  $J= {^{{a_{4}} \ {a_{3}}}_{a_{1}}}$ as an E-pair. We now investigate the non-emptiness of the shift space using the notion of an E-pair.

\begin{Proposition}
Let $X$ be a two dimensional shift of finite type and let the sequence space generated by $M$ and $N$ be non-empty. If for every $M_{ij} \neq 0$ and for every E-pair $" i_{1}"$ of i, $\exists $ an E-pair $"j_{1}"$  of j such that $N_{i_{1}j_{1}} \neq 0 $ 
 then, $X_{G} \neq \phi $.
\end{Proposition}

\begin{proof}
Let $X$ be a shift of finite type such that the sequence spaces generated by $M$ and $N$ are non-empty. Let $M_{ij} \neq 0$ and let $i_1$ be an E-pair of $i$. If there exists an E-pair $"j_{1}"$  of j such that $N_{i_{1}j_{1}} \neq 0$ then the pattern $^{\ \ \  \ {a_{3}} \ ^{ \ a_{5} }_{ \ a_{4}}}_{ \  {a_{1}} \ {a_{2}} }$ can be extended to a $3\times 3$ pattern. As the shift spaces generated by $M$ and $N$ are non-empty, any finite pattern generated by $M$ can be extended to an allowed rectangle of arbitrarily large size and hence can be extended to an element of the shift space. Consequently, the shift space is non-empty and the proof is complete.	
\end{proof}

\begin{Remark}\label{ne2}
The above result establishes the non-emptiness of the shift space using the notion of an E-pair. In particular, the proof establishes that if for every $M_{ij} \neq 0$ and for every E- pair $" i_{1}"$ of i, $\exists $ an E-pair $"j_{1}"$  of j such that $N_{i_{1}j_{1}}\neq 0$ then the shift space is non-empty. It may be noted that the proof ensures the extension of compatible E-pairs into a $3\times 3$ square and hence the condition "for every $M_{ij} \neq 0$ and for every E- pair $" i_{1}"$ of i, $\exists $ an E-pair $"j_{1}"$  of j such that $N_{i_{1}j_{1}} \neq 0 $" is sufficient (but not necessary) to ensure non-emptiness of the shift space. A similar argument proves that if or every $N_{kl} \neq 0$ and for every E- pair $" k_{1}"$ of k, $\exists $ an E-pair $"l_{1}"$ of l such that $M_{k_{1}l_{1}} \neq 0 $ then the shift space is non-empty and hence we get the following corollary.
\end{Remark}

\begin{Cor}
Let $X$ be a two dimensional shift of finite type and let the sequence space generated by $M$ and $N$ be non-empty. If for every $N_{kl} \neq 0$ and for every E- pair $" k_{1}"$ of k, $\exists$ an E-pair $"l_{1}"$ of l such that $M_{k_{1}l_{1}} \neq 0 $ then $X_{G} \neq \phi $.	
\end{Cor}

\begin{proof}
The proof follows from discussions in Remark \ref{ne2}
\end{proof}
	
\begin{Proposition}
	A shift space $X_{G}$ is finite if it follows two conditions: \
	\begin{enumerate}
	\item	M and N are permutation matrices.
	\item Every pattern in $A_{1}$ and $A_{2}$ has unique E-pair.
\end{enumerate}	
\end{Proposition}	

\begin{proof}
Let $X$ be a shift space and let $(1)$ and $(2)$ hold good. Firstly, it may be noted that as $M$ and $N$ are permutation matrices, the shift spaces generated by $M$ and $N$ are finite (union of periodic orbits). Further, as every triangular pattern is uniquely extendable to a $2\times 2$ pattern, every infinite pattern generated by $M$ is uniquely extendable to an element of the shift space. Consequently, the shift space is finite and the proof is complete.
\end{proof}

We now give examples to show that shift space may not be finite if any of the above two conditions are dropped.

\begin{ex}
	Let $X$ be the shift space arising from Figure 6. Then, the adjacency matrices associated with graph  are:
	\begin{figure}[h]
		\includegraphics[height=5.0cm, width=10.0cm]{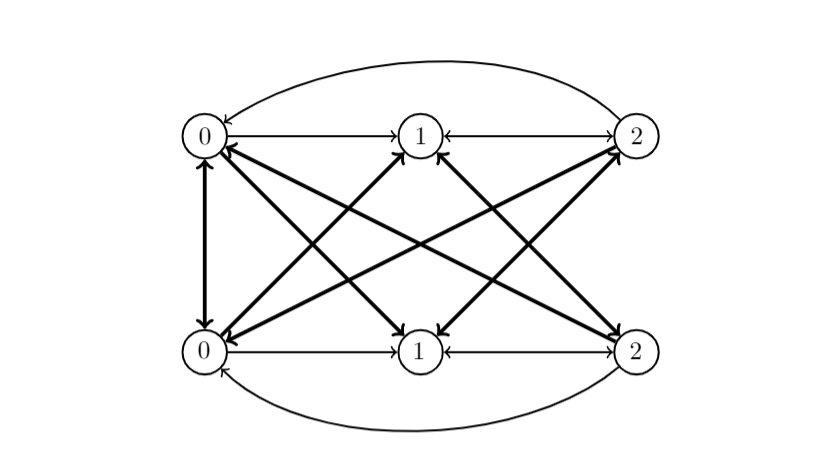}
		\caption{}
	\end{figure}
	
$$\textit{H}= \bordermatrix{ & 0 & 1 & 2 \cr
	0 & 0 & 1 & 0 \cr
	1 & 0 & 0 & 1 \cr
	2 & 1 & 1 & 0 \cr 	
}
\ \ \ \ \ \ \ \ \
\textit{V}= \bordermatrix{ & 0 & 1 & 2 \cr
	0 & 1 & 1 & 0 \cr
	1 & 0 & 0 & 1 \cr
	2 & 1 & 1 & 0 \cr 	
}$$
 Then,
$$\textit{HV}= \bordermatrix{ & 0 & 1 & 2 \cr
	0 & 0 & 0 & 1 \cr
1 & 1 & 1 & 0 \cr
2 & 1 & 1 & 1 \cr
}
\ \ \ \ \ \ \ \ \
\textit{VH}= \bordermatrix{ & 0 & 1 & 2 \cr
	0 & 0 & 1 & 1 \cr
1 & 1 & 1 & 0 \cr
2 & 0 & 1 & 1 \cr 	
}$$
Note that there exists indices $i,j$ such that $(HV)_{ij} \neq 0$ but $(VH)_{ij} =0 $  (and there exists $k, l$ such that $(VH)_{kl} \neq 0 $ but $(HV)_{kl} =0$). Updating the matrices $HV$ and $VH$ we obtain

$$\textit{HV}= \bordermatrix{ & 0 & 1 & 2 \cr
	0 & 0 & 0 & 1 \cr
	1 & 1 & 1 & 0 \cr
	2 & 0 & 1 & 1 \cr 	
}
\ \ \ \ \ \ \ \ \
\textit{VH}= \bordermatrix{ & 0 & 1 & 2 \cr
	0 & 0 & 0 & 1 \cr
1 & 1 & 1 & 0 \cr
2 & 0 & 1 & 1 \cr 		
}$$
Using above matrices, we obtain: \\
$\mathcal{A}_{1} = \{ ^{\ \ \ {2}}_{ \  {0} \ {1} } , ^{\ \ \ {0}}_{ \  {1} \ {2} } , ^{\ \ \ {1}}_{ \  {1} \ {2} } , ^{\ \ \ {1}}_{ \  {2} \ {0} } , ^{\ \ \ {2}}_{ \  {2} \ {1} }  \}$ and
$\mathcal{A}_{2} =\{ {^{{1} \ {2}}_{0}    } , \ {^{{2} \ {0}}_{1}    } , \ {^{{2} \ {1}}_{1}    } , \ {^{{0} \ {1}}_{2}    } , \ {^{{1} \ {2}}_{2}    }  \}$ \\

It can be verified that every element of $\mathcal{A}_{1}$ and $ \mathcal{A}_{2}$ can be extended to $2 \times 2 $ square uniquely and hence condition (2) holds. The matrices M and N are :
$$\textit{M}= \bordermatrix{ & ^{\ \ \ {2}}_{ \  {0} \ {1} } & ^{\ \ \ {0}}_{ \  {1} \ {2} } & ^{\ \ \ {1}}_{ \  {1} \ {2} } &  ^{\ \ \ {1}}_{ \  {2} \ {0} } & ^{\ \ \ {2}}_{ \  {2} \ {1} }  \cr
	^{\ \ \ {2}}_{ \  {0} \ {1} } & 0 & 0 & 0 & 1 & 1  \cr \\
	^{\ \ \ {0}}_{ \  {1} \ {2} } & 1 & 0 & 0 & 0 & 0 \cr \\
	^{\ \ \ {1}}_{ \  {1} \ {2} } & 0 & 1 & 1 & 0 & 0 \cr \\
	^{\ \ \ {1}}_{ \  {2} \ {0} } & 0 & 1 & 1 & 0 & 0\cr
	^{\ \ \ {2}}_{ \  {2} \ {1} } & 0 & 0 & 0 & 1 & 1 \cr
}$$

$$\textit{N}= \bordermatrix{ & {^{{1} \ {2}}_{0}    } &{^{{2} \ {0}}_{1}    } &{^{{2} \ {1}}_{1}    } &   {^{{0} \ {1}}_{2}    } & {^{{1} \ {2}}_{2}    }  \cr
{^{{1} \ {2}}_{0}    } & 0 & 0 & 0 & 1 & 1  \cr \\
{^{{2} \ {0}}_{1}    } & 1 & 0 & 0 & 0 & 0 \cr \\
{^{{2} \ {1}}_{1}    } & 0 & 1 & 1 & 0 & 0 \cr \\
 {^{{0} \ {1}}_{2}    } & 0 & 1 & 1 & 0 & 0\cr
	{^{{1} \ {2}}_{2}    } & 0 & 0 & 0 & 1 & 1 \cr
}$$

It can be seen that $M$ and $N$ are not permutation matrices and the shift space $X$ is not finite. Consequently, the Proposition $8$ does not hold good if $M$ and $N$ are not ensured to be permutation matrices.
\end{ex}

\begin{ex}
	Let $X$ be the shift space arising from graph $G$ in figure-7. Then, the adjacency matrices corresponding to the graph $G$ are:\\
	
\begin{figure}[h]
		\includegraphics[height=5.0cm, width=10.0cm]{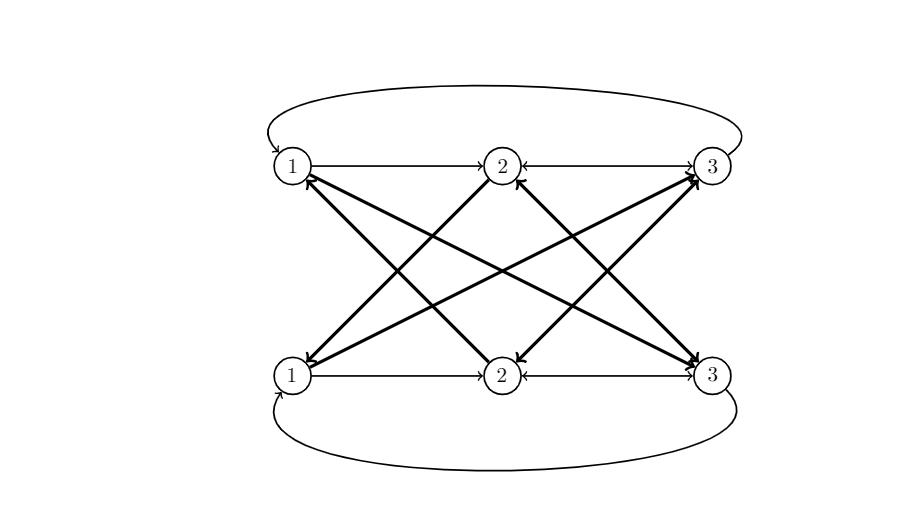}
		\caption{}
	\end{figure}	
	
$$\textit{H}= \bordermatrix{ & 1 & 2 & 3 \cr
	1 & 0 & 1 & 0 \cr
	2 & 0 & 0 & 1 \cr
	3 & 1 & 1 & 0 \cr 	
}
\ \ \ \ \ \ \ \ \
\textit{V}= \bordermatrix{ & 1 & 2 & 3 \cr
	1 & 0 & 0 & 1 \cr
	2 & 1 & 0 & 1 \cr
	3 & 0 & 1 & 0 \cr 	
}$$

Further,

$$\textit{HV}= \bordermatrix{ & 1 & 2 & 3 \cr
	1 & 1 & 0 & 1 \cr
	2 & 0 & 1 & 0 \cr
	3 & 1 & 0 & 2 \cr 	
}
\ \ \ \ \ \ \ \ \
\textit{VH}= \bordermatrix{ & 1 & 2 & 3 \cr
	1 & 1 & 1 & 0 \cr
	2 & 1 & 2 & 0 \cr
	3 & 0 & 0 & 1 \cr 	
}$$
Once again, note that there exists indices $i,j$ such that $(HV)_{ij} \neq 0$ but $(VH)_{ij} =0 $  (and there exists $k, l$ such that $(VH)_{kl} \neq 0 $ but $(HV)_{kl} =0$). Updating the matrices $HV$ and $VH$ we obtain

$$\textit{HV}= \bordermatrix{ & 1 & 2 & 3 \cr
	1 & 1 & 0 & 0 \cr
	2 & 0 & 1 & 0 \cr
	3 & 0 & 0 & 2 \cr 	
}
\ \ \ \ \ \ \ \ \
\textit{VH}= \bordermatrix{ & 1 & 2 & 3 \cr
	1 & 1 & 0 & 0 \cr
	2 & 0 & 2 & 0 \cr
	3 & 0 & 0 & 1 \cr 	
}$$

Consequently, \\
$\mathcal{A}_{1} = \{ ^{\ \ \ {1}}_{ \  {1} \ {2} } , ^{\ \ \ {2}}_{ \  {2} \ {3} } , ^{\ \ \ {3}}_{ \  {3} \ {1} } , ^{\ \ \ {3}}_{ \  {3} \ {2} } \}$ and
$\mathcal{A}_{2} =\{ {^{{3} \ {1}}_{1}    } , \ {^{{3} \ {2}}_{2}    } , \ {^{{1} \ {2}}_{2}    } , \ {^{{2} \ {3}}_{3}    }   \}$ \\

and \\

$$\textit{M}= \bordermatrix{ & ^{\ \ \ {1}}_{ \  {1} \ {2} } & ^{\ \ \ {2}}_{ \  {2} \ {3} } & ^{\ \ \ {3}}_{ \  {3} \ {1} } &  ^{\ \ \ {3}}_{ \  {3} \ {2} }  \cr
	^{\ \ \ {1}}_{ \  {1} \ {2} } & 1 & 0 & 0 & 0   \cr \\
	^{\ \ \ {2}}_{ \  {2} \ {3} } & 0 & 1 & 0 & 0 \cr \\
	^{\ \ \ {3}}_{ \  {3} \ {1} } & 0 & 0 & 1 & 0  \cr \\
	^{\ \ \ {3}}_{ \  {3} \ {2} } & 0 & 0 & 0 & 1 \cr
}$$

$$\textit{N}= \bordermatrix{ & {^{{3} \ {1}}_{1} } & {^{{3} \ {2}}_{2}    } &{^{{1} \ {2}}_{2}    } &  {^{{2} \ {3}}_{3}    }  \cr
	 {^{{3} \ {1}}_{1} } & 1 & 0 & 0 & 0   \cr \\
{^{{3} \ {2}}_{2}    } & 0 & 1 & 0 & 0 \cr \\
	{^{{1} \ {2}}_{2}    } & 0 & 0 & 1 & 0  \cr \\
{^{{2} \ {3}}_{3}    } & 0 & 0 & 0 & 1 \cr
}$$
Clearly, M and N are permutation matrices but not every triangular pattern is getting extended uniquely to $2 \times 2$ pattern. It can be seen that the shift space generated is not finite and hence Proposition $8$ does not hold good if any of the two conditions are dropped.
\end{ex}

\bibliography{xbib}

\end{document}